% To: "Saharon Shelah's Office" <shlhetal@micron.net>
% Subject: Re: F451 -> [HySh:756]
% Date: Mon, 15 Jan 2001 16:25:53 +0200 (EET)
% From: Tapani Hyttinen <thyttine@cc.helsinki.fi>
% Sender: <thyttine@cc.helsinki.fi>
% In-Reply-To: <20001219074224.A2366@dsc235.dsc.unomaha.edu>
% Mime-Version: 1.0
% Content-Description: 
% X-sliced-and-diced-by: 'savemail' 0.4, Jul 2000

\input amssym

\magnification=\magstep1
\hsize=15,4truecm
\vsize=22.5truecm
\advance\voffset by 1truecm
\mathsurround=1pt

\def\chapter#1{\par\bigbreak \centerline{\bf #1}\medskip}

\def\section#1{\par\bigbreak {\bf #1}\nobreak\enspace}

\def\sqr#1#2{{\vcenter{\hrule height.#2pt
      \hbox{\vrule width.#2pt height#1pt \kern#1pt
         \vrule width.#2pt}
       \hrule height.#2pt}}}

\def\k{\kappa}
\def\o{\omega}

\def\d{\delta}

\def\l{\lambda}

\def\a{\alpha}
\def\b{\beta}

\def\g{\gamma}

\def\n{\eta}

%=====================

\def\P{{\cal P}}

%========================

%\def\pm{\buildrel{\scriptscriptstyle +}\over {\scriptscriptstyle -}}
\def\th #1 #2. #3\par\par{\medbreak{\bf#1 #2.
\enspace}{\sl#3\par}\par\medbreak}
\def\rem #1 #2. #3\par{\medbreak{\bf #1 #2.
\enspace}{#3}\par\medbreak}
\def\proof{{\bf Proof}.\enspace}
\def\sqr#1#2{{\vcenter{\hrule height.#2pt
      \hbox{\vrule width.#2pt height#1pt \kern#1pt
         \vrule width.#2pt}
       \hrule height.#2pt}}}
\def\eop{\mathchoice\sqr34\sqr34\sqr{2.1}3\sqr{1.5}3}

%====================================================================%
                                                                     %
% A macro for making references and blocks.                          %
                                                                     %
\newdimen\refindent\newdimen\plusindent                              %
\newdimen\refskip\newdimen\tempindent                                %
\newdimen\extraindent                                                %
                                                                     %
%\refskip has to be determined by the user! Otherwise \parindent is  %
%used, in accordance with \item.                                     %
                                                                     %
\def\ref#1 #2\par{\setbox0=\hbox{#1}\refindent=\wd0                  %
\plusindent=\refskip                                                 %
\extraindent=\refskip                                                %
\advance\extraindent by 30pt                                         %
\advance\plusindent by -\refindent\tempindent=\parindent %           %
\parindent=0pt\par\hangindent\extraindent %                          %
{#1\hskip\plusindent #2}\parindent=\tempindent}                      %
\refskip=\parindent                                                  %
                                                                     %
%====================================================================%

\def\ol{\overline}

\def\empty{\emptyset}

\def\raj{\restriction}

\def\ol{\overline}

\def\nda{\mathrel{\lower0pt\hbox to 3pt{\kern3pt$\not$\hss}\downarrow}}
\def\nDa{\mathrel{\lower0pt\hbox to 3pt{\kern3pt$\not$\hss}\Downarrow}}
\def\nbot{\mathrel{\lower0pt\hbox to 4pt{\kern3pt$\not$\hss}\bot}}
\def\ekom{\mathrel{\lower3pt\hbox to 0pt{\kern3pt$\sim$\hss}\mapsto}}

\def\anR{\mathrel{\lower1pt\hbox to 2pt{\kern3pt$R$\hss}\not}}
\def\anoR{\mathrel{\lower1pt\hbox to 2pt{\kern3pt$\overline{R}$\hss}\not}}

\def\anRm{\mathrel{\lower1pt\hbox to 2pt{\kern3pt$R^{-}$\hss}\not}}

\def\ndda{\mathrel{\lower0pt\hbox to 1pt{\kern3pt$\not$\hss}\downdownarrows}}

\def\warrow{\mathrel{\lower0pt\hbox to 1pt{\kern3pt$^{w}$\hss}\rightarrow}}

\null
\vskip 2truecm
\centerline{\bf FORCING A BOOLEAN ALGEBRA}
\centerline{\bf WITH PREDESIGNED AUTOMORPHISM GROUP}
\vskip 1.5truecm
\centerline{Tapani Hyttinen$^{*}$ and Saharon Shelah$^{\dagger}$}
\vskip 2.5truecm

\chapter{Abstract}

For suitable groups $G$ we will show that one can add a Boolean algebra $B$
by forcing
in such a way that $Aut(B)$ is almost
isomorphic to $G$. In particular, we will
give a positive answer to the following question due to J. Roitman:
Is $\aleph_{\o}$ a possible number of automorphisms of a
rich Boolean algebra?

\vskip 2.5truecm

In [Ro], J. Roitman asked,
is $\aleph_{\o}$ a possible number of automorphisms of a
rich Boolean algebra? A Boolean algebra is rich if the number
of automorphisms is greater than the size of the algebra and
possible means consistent with ZFC (since GCH implies
trivially that the answer is no).
We will answer this question positively. For this we give
a method of adding a Boolean algebra by forcing
in such a way that we have a lot of control
on the automorphism group of the Boolean algebra.
Notice
that by [Mo] Theorem 4.3, one can not hope of
giving a positive answer to Roitmans question
from the assumption that $2^{\o}>\aleph_{\o}$.

Partially, our methods are similar to those in [Ro].

We say that a po-set $P$ is $\k$-Knaster if for all
$p_{i}\in P$, $i<\k$, there is $Z\subseteq\k$ of power $\k$ such that
for all $i,j\in Z$, $p_{i}$ and $p_{j}$ are compatible.

\vskip 1.5truecm

\noindent
$*$ Partially supported by the Academy of Finland,
grant 40734, and the Mittag-Leffler
Institute.

\noindent
$\dagger$ This research was supported by the Israel Science Foundation
founded by the Israel Academy of Sciences and Humanities. Publ. 756

\vfill
\eject

\th Theorem 1. Assume that
$\k =\k^{<\k}$ is an infinite cardinal.
Let $(I,\le )$ be a partial order,
$G_{t}$, $t\in I$, be groups of power $\le\k$ and
$\pi_{s,t}:G_{t}\rightarrow G_{s}$, $s <t\in I$, be homomorphisms
such that

(a) $I$ is $\k^{+}$-directed,
any pair of elements of $I$
has the least upper bound and for all $t\in I$,
$\vert\{ s\in I\vert\ s\le t\}\vert\le\k$,

(b) for $r<s<t\in I$, $\pi_{r,t}=\pi_{r,s}\circ\pi_{s,t}$.

\noindent
Then there is a $\k$-closed, $\k^{+}$-Knaster
po-set $P$ such that
in $V^{P}$ the following holds:
There is a Boolean algebra $B$ such that

(i) $B$ is atomic with $\vert I\vert +\k$ atoms,

(ii) $\vert B\vert =\k^{+}+\vert I\vert$,

(iii) if $G^{*}$ is the inverse limit of $(G_{t},\pi_{s,t}\vert s<t\in I)$,
then $G^{*}$ can be embedded into $Aut(B)$,
$\vert Aut(B)\vert\le\k^{+}\cdot\vert I\vert\cdot\vert G^{*}\vert$
and if for all $t\in I$, $G_{t}$ is the one element group, then
$\vert Aut(B)\vert =\k\cdot\vert I\vert$.

Note: We will show that $Aut(B)$ is not very far from
$G^{*}$. In particular, we will show that
for all $f\in Aut(B)$ there is $g\in G^{*}$ such that
excluding a small error,
$f$ is the same as the image of
$g$ under the embedding.

\proof For $t\in I$, let $X_{t}=\{ t\}\times G_{t}\times\k$,
$X_{\le t}=\cup_{s\le t}X_{s}$ and $X=\cup_{t\in I}X_{t}$.
We define $P$ to be the set of all tuples $p=(U,T,V,Y)=
(U^{p},T^{p},V^{p},Y^{p})$ such that

(1) $U\subseteq X$ is of power $<\k$,

(2) $T\subseteq I$ is of power $<\k$,

(3) $V=(v_{t}\vert\ t\in T)=(v^{p}_{t}\vert\ t\in T)$, where each
$v_{t}\subseteq\k^{+}$ is of power $<\k$,

(4) $Y=(y_{t,\a}\vert\ t\in T,\ \a\in v_{t})=
(y^{p}_{t,\a}\vert\ t\in T,\ \a\in v_{t})$, where
$y_{t,\a}\subseteq U\cap X_{\le t}$.

\noindent
We order $P$ so that $p\le q$ ($q$ is stronger than $p$) if

(5) $U^{p}\subseteq U^{q}$,

(6) $T^{p}\subseteq T^{q}$,

(7) $v^{p}_{t}\subseteq v^{q}_{t}$ for $t\in T^{p}$,

(8) $y^{p}_{t,\a}=y^{q}_{t,\a}\cap U^{p}$ for
$t\in T^{p}$ and $\a\in v^{p}_{t}$.

Clearly $P$ is $\k$-closed and $\k^{<\k}=\k$ implies
by the usual $\Delta$-lemma argument that
$P$ is $\k^{+}$-Knaster.

Let $H$ be $P$-generic over $V$. We let
$Y_{t,\a}=\cup\{ y^{p}_{t,\a}\vert\ p\in H\}$.
Below, we will not distinguish
these and other objects from their names, it will be clear from the
context which we mean.

Let $G^{*}$ be as in the theorem i.e. $G^{*}$ consists
of all $g:I\rightarrow\cup_{t\in I}G_{t}$ such that
for all $s<t\in I$, $\pi_{s,t}(g(t))=g(s)$.
If $g\in G^{*}$
and $x=(t,z,\a )\in X_{t}$, then we write
also $g(x)$ for $(t,g(t)z,\a )$
and define $Y_{t,\a ,g}=\{ g(x)\vert\ x\in Y_{t,\a}\}$.

By a Boolean term $\d (Y_{0},...,Y_{n-1})$ taken inside
$Y$ we mean a term of the form
$\cap_{i<n}Z_{i}$, where each $Z_{i}$ is either
$Y\cap Y_{i}$ or $Y-Y_{i}$. If $Z_{i}=Y\cap Y_{i}$ we say that
$Y_{i}$ appears in $\d$ positively. A general Boolean
term is a finite union of Boolean terms.

\th 1.1 Claim. Let $\d (Y_{0},...,Y_{n-1})$
be a Boolean term taken inside $X_{\le t}$, $t\in I$ and
for $i<n$, $t_{i}\in I$, $\a_{i}<\k^{+}$ and $g_{i}\in G^{*}$.
Assume that
there is $t'\in I$ such that

(a) $t'\le t$,

(b) if $Y_{i}$ appears in $\d$
positively, then $t'\le t_{i}$,

(c) if $Y_{i}$ appears in $\d$ positively and $Y_{j}$
does not appear positively in $\d$, then
$(t_{i},\a_{i},g_{i}(t'))\ne (t_{j},\a_{j},g_{j}(t'))$.

\noindent
Then $\d (Y_{t_{0},\a_{0},g_{0}},...,Y_{t_{n-1},\a_{n-1},g_{n-1}})
\cap X_{\le t'}$
has power $\k$.
On the other hand, if such a $t'$ does not exist, then
$\d (Y_{t_{0},\a_{0},g_{0}},...,Y_{t_{n-1},\a_{n-1},g_{n-1}})=\empty$.

\proof An easy density argument using the fact that
the actions of two distinct
elements of a group differ on each element of the group.
$\eop$ Claim 1.1

Now we can define the Boolean algebra $B$. We let $B$
be the subalgebra of $\P (X)$ generated by
$[X]^{<\k}\cup\{ Y_{t,\a ,g}\vert\ t\in I,\ \a <\k^{+}, g\in G^{*}\}$,
where $[X]^{<\k}$ is the set of all subsets of $X$ of power $<\k$.

\th 1.2 Claim. Items (i) and (ii) from Theorem 1 hold
and $G^{*}$ can be embedded into $Aut(B)$.

\proof (i) is clear.
For (ii), notice that $Y_{t,\a ,g}$ is determined by
$Y_{t,\a}$ and $g(t)$.
Finally, it is easy to check
that the function which takes $g\in G^{*}$ to
$\ol g:A\mapsto\{ g(z)\vert\ z\in A\}$ is the required
embedding of $G^{*}$ into $Aut(B)$. $\eop$ Claim 1.2.

So we are 'left' to prove the second half of (iii) from Theorem 1:
For all $0<m<\o$ and $t\in I$, let $\d^{t,m}_{i}(Y_{0},...,Y_{m-1})$,
$i<2^{m}$, list all Boolean terms in $m$ variables taken
inside $X_{\le t}$
so that
$\d^{t,m}_{0}(Y_{0},...,Y_{m-1})=\cap_{k<m}(X_{\le t}-Y_{k})$
and for $m=0$, we let $\d^{t,m}_{0}=X_{\le t}$.

Let $h\in Aut(B)$. By the definition of $B$,
$h$ is determined by $h\raj\{\{ x\}\vert\ x\in X\}$.
In order to simplify the notation, for $x\in X$, we write
$h(x)$ for the unique $z\in X$ such that $h(\{ x\} )=\{ z\}$.
We will frequently consider such a $h$ as if it is a permutation of $X$.

\th 1.3 Claim. In $V^{P}$, if $h\in Aut(B)$ and $t\in I$,
then there are $u\subseteq X_{\le t}$, $m<\o$,
for $i<m$, $t_{i}\in I$, $\a_{i}<\k^{+}$ and
$g_{i}\in G^{*}$ and for $j<2^{m}$, $f_{j}\in G^{*}$
such that

(i) $\vert u\vert <\k$,

(ii) either for all $j<2^{m}$ and
$x\in\d^{t,m}_{j}(Y_{t_{0},\a_{0},g_{0}},...,Y_{t_{m-1},\a_{m-1},g_{m-1}})-u$,
$h(x)=f_{j}(x)$ or for all $j<2^{m}$ and
$x\in\d^{t,m}_{j}(Y_{t_{0},\a_{0},g_{0}},...,Y_{t_{m-1},\a_{m-1},g_{m-1}})-u$,
$h^{-1}(x)=f_{j}(x)$.

\proof For a contradiction, assume that $p$ forces
that $h$ and $t$ are a counter example.
For all $x\in X_{\le t}$,
let $(p_{x,i})_{i<\k}$ be a maximal antichain
such that for all $i<\k$, $p_{x,i}\ge p$ or
$p_{x,i}$ is incompatible with $p$ and
$p_{x,i}$ decides both $h(x)$ and $h^{-1}(x)$.
Choose $\a^{*}<\k^{+}$ so that
$\a^{*}\not\in\cup\{ v^{p_{x,i}}_{t}\vert\ x\in X_{\le t},\ i<\k\}$.
Choose $q\ge p$ so that
it decides both $h(Y_{t,\a^{*}})$ and $h^{-1}(Y_{t,\a^{*}})$ i.e.
it forces $h(Y_{t,\a^{*}})$
to be, say, $(u_{0}\cup\d (Y_{t'_{0},\a'_{0},g'_{0}},...,
Y_{t'_{k-1},\a'_{k-1},g'_{k-1}}))-u_{1}$ and
$h^{-1}(Y_{t,\a^{*}})$
to be $(u^{-}_{0}\cup\d^{-}(Y_{t'_{0},\a'_{0},g'_{0}},...,
Y_{t'_{k-1},\a'_{k-1},g'_{k-1}}))-u^{-}_{1}$,
where $u_{0},u_{1},u^{-}_{1},u^{-}_{2}\in [X]^{<\k}$ and
$\d$ and $\d^{-}$ are general Boolean terms taken inside $X$.
(Notice that $Y_{t,\a^{*}}=Y_{t,\a^{*},1}$.)
In particular, $q$ is chosen so that it decides values for
$g'_{i}(t'_{i})$ and $g'_{i}(t)$, $i<k$.
Clearly we may assume that ($q$ forces that)
in $((t'_{0},\a'_{0},g'_{0}(t'_{0})),...,
(t'_{k-1},\a'_{k-1},g'_{k-1}(t'_{k-1})))$ there is no repetition,
and for some $k^{*}\le k$,
$\{ i<k\vert\ (t'_{i},\a'_{i})=(t,\a^{*})\} =\k^{*}$.
We write $\ol Y_{<k^{*}}$ for $(Y_{t'_{0},\a'_{0},g'_{0}},...,
Y_{t'_{k^{*}-1},\a'_{k^{*}-1},g'_{k^{*}-1}})$
and $\ol Y_{\ge k^{*}}$ for $(Y_{t'_{k^{*}},\a'_{k^{*}},g'_{k^{*}}},...,
Y_{t'_{k-1},\a'_{k-1},g'_{k-1}})$.

If we can find $x\in X_{\le t}-U^{q}$ and $q'\ge q$ such that
$q'$ forces that
for all $i<k^{*}$, $g'_{i}(x)\ne h(x)$ and that $h(x)\not\in u_{0}\cup u_{1}$,
we get an easy
contradiction.
The contradiction is got the same way as we get a contradiction
after the proof of Subclaim 1.3.1 (we can find $q'\ge q$ so that
it forces
$h(x)\in h(Y_{t,\a^{*}})\ \leftrightarrow\ h(x)\not\in\d (\ol Y)$.)
Similarly, if we can find $x\in X_{\le t}-U^{q}$ and $q'\ge q$ such that
$q'$ forces that
for all $i<k^{*}$, $g'_{i}(x)\ne h^{-1}(x)$ and
that $h^{-1}(x)\not\in u^{-}_{0}\cup u^{-}_{1}$,
we get an easy
contradiction.
So we assume the following:

(*) For all $x\in X_{\le t}-U^{q}$, $q$ forces that
$h(x)\in u_{0}\cup u_{1}$ or
for some $i<k^{*}$, $h(x)=g'_{i}(x)$ and that
$h^{-1}(x)\in u^{-}_{0}\cup u^{-}_{1}$ or
for some $i<k^{*}$, $h^{-1}(x)=g'_{i}(x)$.

Notice that (*) implies that
for all $t'\le t$, $q$ forces that
excluding a small error,
$h$ is a permutation of $X_{t'}$.

By (*) above and the choice of $p$, $h$ and $t$,
and by making $q$ stronger if necessary,
we may assume that there are
Boolean terms $\d'$ and $\d''$ taken inside $X_{\le t}$
such that 

($\a$) $q$ forces that for all $g\in G^{*}$ there are $\k$ many
$x\in X_{\le t}$
such that $g(x)\ne h(x)\in\d'(\ol Y_{<k^{*}})\cap
\d''(\ol Y_{\ge k^{*}})$,

($\b$) $q$ forces that
$\d'(\ol Y_{<k^{*}})\cap
\d''(\ol Y_{\ge k^{*}})\ne\empty$,

($\g$) $\empty$ forces that
$\d'(\ol Y_{<k^{*}})\cap
\d''(\ol Y_{\le k^{*}})\subseteq\d (\ol Y)$.

\noindent
Instead of item ($\g$), we may have
$\d'(\ol Y_{<k^{*}})\cap
\d''(\ol Y_{\le k^{*}})\subseteq X_{\le t}-\d (\ol Y)$. This
is a symmetric case and we skip it.

\th 1.3.1 Subclaim. There are $x_{1},x_{2}\in X_{\le t}$,
$d_{1},d_{2}\in X_{\le t}$,
$q^{*}\ge q$ and non-empty $W_{1},W_{2}\subseteq k^{*}$
such that

(1) $W_{1}\cap W_{2}=\empty$,

(2) $q^{*}$ forces that for all $n\in\{ 1,2\}$, $h(x_{n})=d_{n}$,

(3) $q^{*}$ forces that for all $n\in\{ 1,2\}$,
$d_{n}\in \d''(\ol Y_{\ge k^{*}})-(u_{0}\cup u_{1})$,

(4) $q^{*}$ forces that for all $n\in\{ 1,2\}$ and $i<k^{*}$,
$d_{n}=g'_{i}(x_{n})$ iff $i\in W_{n}$,

(5) if $d_{1}\in X_{v_{1}}$ and $d_{2}\in X_{v_{2}}$,
then $v_{1}\le v_{2}$,

(6) for all $n\in\{ 1,2\}$ and $i<k^{*}$,
$(g'_{i})^{-1}(d_{n})\not\in U^{q}$ and
$\{ (g'_{i})^{-1}(d_{1})\vert\ i<k^{*}\}\cap
\{ (g'_{i})^{-1}(d_{2})\vert\ i<k^{*}\} =\empty$.

\proof Let $I'\subseteq I$ be the set of those $t'\le t$
such that $X_{t'}\cap\d''(\ol Y_{\le k^{*}})\ne\empty$.
Notice that by Lemma 1.1, $q$ decides $I'$ and
$I'$ is closed under the least upper bounds.
By ($\b$) above, $I'\ne\empty$.
For $t'\in I'$, we write $I'_{\ge t'}$ for
$\{ v\in I'\vert\ v\ge t'\}$.
Then for all
$t^{*}\in I'$ there are $W^{*}\subseteq k^{*}$ and $q'\ge q$,
such that $q'$ forces that excluding $<\k$ many $x$ if
$h(x)\in \d''(\ol Y_{\ge k^{*}})-(u_{0}\cup u_{1})$ and
$x\in X_{v}$
for some $v\in I'_{\ge t'}$,
then
$h(x)=g'_{i}(x)$ for some $i\in W^{*}$.
We choose these so that $\vert W^{*}\vert$ is minimal.
There are two cases:

(I) $\vert W^{*}\vert >1$: Now it is easy to choose
$v_{1},v_{2}\in I'_{\ge t^{*}}$, $x_{1}\in X_{v_{1}}$,
$x_{2}\in X_{v_{2}}$, $d_{1},d_{2}\in X_{\le t}$,
$i\in W^{*}$
and $q^{*}\ge q'$ so that
$v_{2}\ge v_{1}$, for $n\in\{ 1,2\}$, $q^{*}$ forces that
$h(x_{n})=d_{n}$,
$h(x_{n})\in \d''(\ol Y_{\ge k^{*}})-(u_{0}\cup u_{1})$,
$h(x_{1})\ne g'_{i}(x_{1})$ and $h(x_{2})=g'_{i}(x_{2})$.
Then we let $W_{1}=\{ j<k^{*}\vert g'_{j}(x_{1})=d_{1}\}$
and $W_{2}=\{ j<k^{*}\vert g'_{j}(x_{2})=d_{2}\}$.
By the choice of $W^{*}$, these can
be chosen so that in addition (6) holds.
These are as wanted ($v_{2}\ge v_{1}$ and $h(x_{1})\ne g'_{i}(x_{1})$
imply (1)).

(II) $\vert W^{*}\vert =1$: Let $i^{*}$ be the only element
of $W^{*}$.
By ($\a$) above, we can find $q''\ge q'$, $x_{1}\in X_{\le t}-U^{q}$
and $d_{1}\in X_{\le t}$,
so that $q''$ forces that $h(x_{1})=d_{1}$,
$h(x_{1})\in \d''(\ol Y_{\ge k^{*}})-(u_{0}\cup u_{1})$
and $h(x_{1})\ne g'_{i^{*}}(x_{1})$. Let $v_{1}\in I$ be such that
$x_{1}\in X_{v_{1}}$ and let $v_{2}$ be the least upper bound
of $v_{1}$ and $t^{*}$. By the assumption (*) above,
$d_{1}\in X_{v_{1}}$ and so
$v_{2}\in I'_{\ge t^{*}}$.
Then by Claim 1.1 and the assumption
that $\vert W^{*}\vert =1$, there are
$x_{2}\in X_{v_{2}}$, $d_{2}\in X_{\le t}$ and $q^{*}\ge q''$ such that
$q^{*}$ forces that $h(x_{2})=d_{2}=g'_{i^{*}}(x_{2})$ and
$h(x_{2})\in \d''(\ol Y_{\ge k^{*}})-(u_{0}\cup u_{1})$.
Then we can let $W_{1}=\{ j<k^{*}\vert g'_{j}(x_{1})=d_{1}\}$
and $W_{2}=\{ j<k^{*}\vert g'_{j}(x_{2})=d_{2}\}$.
By ($\a$) and the choice of $W^{*}$, these can
be chosen so that in addition (6) holds.
These are as wanted. $\eop$ Subclaim 1.3.1.

By the choice of the conditions $p_{x,i}$, we can find
$r'\ge q^{*}$ and $i,i'<\k$ so that
$p_{x_{1},i},p_{x_{2},i'}\le r'$.

We define a condition $r$ as follows:

(a) $U^{r}=U^{r'}\cup\{ (g'_{i})^{-1}(d_{n})\vert\ i<k,\ n\in\{ 1,2\}\}$,

(b) $T^{r}=T^{r'}\cup\{ t\}$,

(c) $v^{r}_{s}=v^{r'}_{s}\cup\{\a^{*}\}$,

(d) if $(s,\a )\ne (t,\a^{*})$, then $y^{r}_{s,\a}=y^{r'}_{s,\a}$,

(e) $y^{r}_{t,\a^{*}}=y^{q}_{t,a^{*}}\cup
\{  (g'_{i})^{-1}(d_{n})\vert\ i\in W_{1},\ n\in\{ 1,2\}\}$.

\noindent
Then $r\ge q,p_{x_{1},i},p_{x_{2},i'}$, in particular
$r$ forces that $h(x_{1})=d_{1}$ and $h(x_{2})=d_{2}$.
Let $\d^{*}(Y_{0},...,Y_{k^{*}-1})$ be a boolean term
taken inside $X_{\le t}$ so that $Y_{i}$ appears in $\d^{*}$
positively iff $i\in W_{1}$. Then
using Subclaim 1.3.1,
it is easy to check
that $r$ forces that
$d_{1}\in\d^{*}(\ol Y_{<k^{*}})\cap\d''(\ol Y_{\ge k^{*}})-(u_{0}\cup u_{1})$
and $d_{1}\in h(Y_{t,\a^{*}})$. Since $r\ge q$,
this implies that $r$
forces that
$\d^{*}(\ol Y_{<k^{*}})\cap\d''(\ol Y_{\ge k^{*}})-u_{1}
\subseteq
\d (\ol Y)-u_{1}\subseteq
h(Y_{t,\a^{*}})$. On the other hand, using Subclaim 1.3.1 (1), (4) and (5),
$r$ forces that
$d_{2}\in\d^{*}(\ol Y_{<k^{*}})\cap\d''(\ol Y_{\ge k^{*}})-(u_{0}\cup u_{1})$
and $d_{2}\not\in h(Y_{t,\a^{*}})$, a contradiction.
$\eop$ Claim 1.3.

\relax From now on we work in $V^{P}$.

\th 1.4 Claim. Suppose $h\in Aut(B)$.
There are $h^{*}$,
$s^{*}\in I$ $u^{s^{*}}\subseteq X_{\le s^{*}}$, $m^{*}<\o$,
for $i<m^{*}$, $t^{s^{*}}_{i}\in I$,
$\a^{s^{*}}_{i}<\k^{+}$, $g^{s^{*}}_{i}\in G^{*}$,
for $0<j<2^{m^{*}}$, $f^{s^{*}}_{j}\in G^{*}$ and
$f^{*}\in G^{*}$
such that if we
write $\ol Y=(Y_{t^{s^{*}}_{0},\a^{s^{*}}_{0},g^{s^{*}}_{0}},...,
Y_{t^{s^{*}}_{m^{s^{*}}-1},\a^{s^{*}}_{m^{s^{*}}-1},g^{s^{*}}_{m^{s^{*}}-1}})$,
then the following holds:

(i) $h^{*}=h$ or $h^{*}=h^{-1}$,

(ii) for all $0<j<2^{m^{*}}$ and $x\in\d^{s^{*},m^{*}}_{j}(\ol Y)-u^{*}$,
$h^{*}(x)=f^{s^{*}}_{j}(x)$,

(iii) for all
$x\in (X-X_{\le s^{*}})\cup (\d^{s^{*},m^{*}}_{0}(\ol Y)-u^{*})$,
$h^{*}(x)=f^{*}(x)$,

(iv) for all $i<m^{*}$, $t^{s^{*}}_{i}\le s^{*}$,

(v) $\vert u^{*}\vert <\k$.

\proof For every $s\in I$, choose $u^{s}$, $m^{s}$,
$t^{s}_{i}$, $\a^{s}_{i}$, $g^{s}_{i}$ and $f^{s}_{j}$,
$i<m^{s}$ and $j<2^{m^{s}}$, as in Claim 1.3.
We let $I^{*}$ to be the set of all
$s\in I$ such that
for all $j<2^{m^{s}}$ and
$x\in\d^{s,m^{s}}_{j}
(Y_{t^{s}_{0},\a^{s}_{0},g^{s}_{0}},...,
Y_{t^{s}_{m-1},\a^{s}_{m-1},g^{s}_{m-1}})-u^{s}$,
$h(x)=f^{s}_{j}(x)$.
We may assume that $I^{*}$ is unbounded
in $I$ (the case when $I-I^{*}$ is unbounded is similar) and we let
$h^{*}=h$.
We may also assume that these are chosen so that in addition
$m^{s}$ and $u^{s}$ are minimal (in this order).

We write $E(s)$, $s\in I^{*}$, for the set of all
$(t,\a )$ such that for some $i<m^{s}$, $t=t^{s}_{i}$ and
$\a =\a^{s}_{i}$.

\th 1.4.1 Subclaim. If $r,s\in I^{*}$ and $r<s$, then

(a) $E(r)\subseteq E(s)$,

(b) for all $j<2^{m^{r}}$, $k<2^{m^{s}}$ and $t\le r$,
if
$$\d^{r,m^{r}}_{j}(Y_{t^{r}_{0},\a^{r}_{0},g^{r}_{0}},...,
Y_{t^{r}_{m^{r}-1},\a^{r}_{m^{r}-1},g^{r}_{m^{r}-1}})\cap$$
$$\d^{s,m^{s}}_{k}(Y_{t^{s}_{0},\a^{s}_{0},g^{s}_{0}},...,
Y_{t^{s}_{m^{s}-1},\a^{s}_{m^{s}-1},g^{s}_{m^{s}-1}})
\cap X_{t}\ne\empty,$$
then $f^{r}_{j}(t)=f^{s}_{k}(t)$.

\proof We will prove (a), (b) is immediate by Lemma 1.1.
For a contradiction, assume that

(1) $(t^{r}_{0},\a^{r}_{0})
\not\in E(s)$.

We will write $\ol Y^{r}=(Y_{t^{r}_{0},\a^{r}_{0},g^{r}_{0}},...,
Y_{t^{r}_{m^{r}-1},\a^{r}_{m^{r}-1},g^{r}_{m^{r}-1}})$ and
$\ol Y^{r}_{*}=(Y_{t^{r}_{1},\a^{r}_{1},g^{r}_{1}},...,
$ $Y_{t^{r}_{m^{r}-1},\a^{r}_{m^{r}-1},g^{r}_{m^{r}-1}})$ and
similarly for $s$ instead of $r$.

We will show that,
we can find $a,b<2^{m^{r}}$, $t\le r$ and a Boolean term $\d$ taken
inside $X_{\le r}$ such that

(2) $f^{r}_{a}(t)\ne f^{r}_{b}(t)$,

(3) $\d^{r,m^{r}}_{a}(\ol Y^{r})=
(X_{\le r}\cap Y_{t^{r}_{0},\a^{r}_{0},g^{r}_{0}})
\cap\d (\ol Y^{r}_{*})$ and
$X_{t}\cap\d^{r,m^{r}}_{a}(\ol Y^{r})\ne\empty$,

(4) $\d^{r,m^{r}}_{b}(\ol Y^{r})=(X_{\le r}-Y_{t^{r}_{0},\a^{r}_{0},g^{r}_{0}})
\cap\d (\ol Y^{r}_{*})$ and $X_{t}\cap\d^{r,m^{r}}_{b}(\ol Y^{r})\ne\empty$.

\noindent
If not, then for all $i<2^{m^{r}-1}$, we define $f_{i}\in G^{*}$
so that for all $x\in\d^{r,m^{r}-1}_{i}(Y^{r}_{*})-u^{r}$,
$h(x)=f_{i}(x)$, this contradicts the minimality of $m^{r}$.

Let $J$ be the set of all $t\le t^{r}_{0}$ such that
$X_{t}\cap\d^{r,m^{r}-1}_{i}(\ol Y^{r}_{*})\ne\empty$.
Notice that by Lemma 1.1,
$J$ is closed under the least upper bounds.
Let $c,d<2^{m^{r}}$ be such that
$\d^{r,m^{r}}_{c}(\ol Y^{r})=
(X_{\le r}\cap Y_{t^{r}_{0},\a^{r}_{0},g^{r}_{0}})
\cap\d^{r,m^{r}-1}_{i}(\ol Y^{r}_{*})$
and
$\d^{r,m^{r}}_{d}(\ol Y^{r})=(X_{\le r}-Y_{t^{r}_{0},\a^{r}_{0},g^{r}_{0}})
\cap\d^{r,m^{r}-1}_{i}(\ol Y^{r}_{*})$.
Let $J_{c}$ be the set of all $t\in J$ such that
$X_{t}\cap\d^{r,m^{r}}_{c}(\ol Y^{r})\ne\empty$.
$J_{d}$ is defined similarly. Then either
$J_{c}$ or $J_{d}$ is cofinal in $J$.
If $J_{d}$ is cofinal, then we let $f_{i}=f^{r}_{d}$
and otherwise we let $f_{i}=f^{r}_{c}$.

We show that $f_{i}$ is as wanted. We may assume that $J_{d}$
is not cofinal, the other case is easy.
Then since $J$ is closed under the least upper bounds,
$J-J_{d}$ is cofinal. By Lemma 1.1, this means that there is
$0<k<m^{r}$ such that for all $t\in J-J_{d}$,
$(t^{r}_{k},\a^{r}_{k},g^{r}_{k}(t))=(t^{r}_{0},\a^{r}_{0},g^{r}_{0}(t))$
and $Y_{t^{r}_{k},\a^{r}_{k},g^{r}_{k}}$ appears
positively in $\d^{r,m^{r}-1}_{i}(\ol Y^{r}_{*})$.
But then $\d^{r,m^{r}-1}_{i}(\ol Y^{r}_{*})\subseteq
Y_{t^{r}_{0},\a^{r}_{0},g^{r}_{0}}$ and so $f_{i}$ is as wanted.

We have shown the existence of
$a,b,t$ and $\d$.

By Claim 1.1 and the fact that $s\ge r$, there is $c<2^{m^{s}}$
such that $\d^{s,m^{s}}_{c}(\ol Y^{s})\cap\d^{r,m^{r}}_{a}(\ol Y^{r})
\cap X_{t}\not\subseteq u^{s}$. But then, by (1) above and Claim 1.1,
$\d^{s,m^{s}}_{c}(\ol Y^{s})\cap\d^{r,m^{r}}_{b}(\ol Y^{r})
\cap X_{\le t}\not\subseteq u^{s}$. This contradicts (2).
$\eop$ Subclaim 1.4.1.

Since $I$ is $\k^{+}$-directed, also
$I^{*}$ is $\k^{+}$-directed and so Subclaim 1.4.1 implies that there
is $s^{*}$ such that
$E(s^{*})$ is maximal and for all $i<m^{*}=m^{s^{*}}$,
$t^{s^{*}}_{i}\le s^{*}$. By the choice of $u^{s}$, $s\in I^{*}$,
it is easy to see that if $r,s\in I^{*}$ and $s^{*}\le r\le s$,
then $u^{r}\subseteq u^{s}$. So
$s^{*}$ can be chosen so
that, in addition,
for all $s^{*}\le r\in I^{*}$, $u^{r}=u^{s^{*}}$.
By Claim 1.1 and Subclaim 1.4.1,
there is $f^{*}\in G^{*}$ such that for all
$s\in I^{*}$, $f^{s}_{0}(s)=f^{*}(s)$.
These are as wanted.
$\eop$ Claim 1.4.

Now it follows immediately from Claim 1.4 that
$\vert Aut(B)\vert\le\k^{+}\cdot\vert I\vert\cdot\vert G^{*}\vert$
and if for all $t\in I$, $G_{t}$ is the one element group, then
$\vert Aut(B)\vert\le\k\cdot\vert I\vert$.
Finally, for $\vert Aut(B)\vert\ge\k\cdot\vert I\vert$,
notice that every permutation of $X$ which is identity on
every element except on $<\k$ many, can be lifted to an
automorphism of $B$.
$\eop$

\th 2 Remark. Notice that in the proof of Theorem 1, one can add
$\l>\k^{+}$ subsets for the sets $X_{\le t}$ making the Boolean
algebra larger,
and the proof still
works giving $\vert G^{*}\vert\le\vert Aut(B)\vert\le\l\cdot\vert I\vert
\cdot\vert G^{*}\vert$.

Our aim was to use
Theorem 1 to solve questions on the number of
automorphisms of Boolean algebras.
We start by proving first a result
from [Ro].

\th 3 Conclusion. ([Ro]) Assume that $\k^{<\k}=\k$. Then there is
a $\k$-closed, $\k^{+}$-Knaster po-set $P$ of power $\k^{+}$ such that in
$V^{P}$ there is an atomic Boolean algebra $B$ of power $\k^{+}$
with $\vert Aut(B)\vert =\k$.

\proof Just choose $I=\{ 0\}$ and $G_{0}$ to be
the one element group and apply Theorem 1
(with these choices of $I$ and $G_{0}$, our proof
is essentially the same as the related proof in [Ro]). $\eop$

We say that $T$ is a $\k$-tree if it is a tree  of power $\k$
and of height $\k$.
We write $T_{\a}$ for $\{ t\in T\vert\ otp(\{ s\in T\vert\ s<t\} )=\a\}$
and say that a $\k$-tree $T$ is narrow if for all
$\a<\k$, $\vert T_{\a}\vert <\k$.
We say that $b$ is a branch in a tree $T$,
if $b$ is a linearly ordered downwards closed subtree of $T$
and the height of $b$ is the same as that of $T$.

\th 4 Corollary. Assume that $\k^{<\k}=\k$ and that there is
a narrow $\k^{+}$-tree with $\l >0$ branches. Then there is
a $\k$-closed, $\k^{+}$-Knaster po-set $P$ of power $\k^{+}$
with the following property: in
$V^{P}$ there is an atomic Boolean algebra $B$ of power $\k^{+}$
such that $\l\le\vert Aut(B)\vert\le\k^{+}\cdot\l$.

\proof We apply Theorem 1:
Clearly we may assume that every element of $T$
belongs to some branch. Then we
let $I=\k^{+}$ and for each $\a\in I$, let
$G_{\a}$ be the free group generated by
$T_{\a}$.
Finally, let $\pi_{\a ,\b}$ be the natural projection
along the branches of $T$.
Then by Theorem 1, it is enough to show the following:

\th 4.1 Claim. The cardinality of the
inverse limit of the system is $\o\cdot\l$ (in $V^{P}$).

\proof It is enough to show that the forcing does not add
branches to $T$ (since $cf(\k^{+})>\o$, the inverse limit
is the free group generated by the 'branches' of $T$).
For a contradiction assume that
$p$ forces that $b$ is a new branch in $T$.
For all $\a <\k^{+}$, choose $p_{\a}\ge p$ and $t_{\a}\in T_{\a}$
so that $p$ forces that $t_{\a}\in b$.
Since $P$ is
$\k^{+}$-Knaster, there is $Z\subseteq\k^{+}$ of power
$\k^{+}$ such that for all $\a ,\b\in Z$,
$p_{\a}$ and $p_{\b}$ are compatible (and so if $\a <\b$, then
$t_{\a}<t_{\b}$). Then
there are no $\a\in Z$ so that $p_{\a}$ forces
that there are $\k^{+}$ many $\b\in Z$ such that $p_{\b}$
is in the generic set.
So for every $\g <\k^{+}$ we can find $i_{\g}<\k^{+}$ and
$q_{\g}\ge p_{i_{\g}}$ such that $i_{\g}\in Z$,
if $\g <\g'$, then $i_{\g}<i_{\g'}$ and
for all $j\ge i_{\g +1}$, if $j\in Z$, then
$q_{\g}$ forces that $p_{j}$ does not belong to the
generic set. But then $\{ q_{\g}\vert\ \g <\k^{+}\}$
is an antichain, a contradiction. $\eop$ Claim 4.1.

$\eop$

\th 5 Conclusion. Con(ZFC) implies the consistency of
ZFC together with the following: 
there is an atomic Boolean algebra $B$
of power $\aleph_{1}$ such that $\vert Aut(B)\vert =\aleph_{\o}$
(and if one wants $2^{\aleph_{0}}>\aleph_{\o}$).

\proof By Corollary 4 it is enough to show that Con(ZFC) implies
the consistency of ZFC together with the following:
($2^{\aleph_{0}}>\aleph_{\o}$, if wanted,
and) there is a narrow $\o_{1}$-tree $T$ with $\aleph_{\o}$ branches.
This is standard but let us sketch the proof: Assume CH.
Let $P$ consist of four tuples $p=(\a^{p},T^{p},u^{p},
\{\n^{p}_{i}\vert\ i\in u^{p}\} )$
such that $\a^{p}<\o_{1}$, $T^{p}$ is a
downwards closed countable subtree of $\o^{\le\a}$,
$u^{p}$ is a countable
non-empty subset of $\aleph_{\o}$ and for all
$i\in u^{p}$, $\n^{p}_{i}$ is a branch in $T^{p}$.
The order is the obvious one: $p<q$ if
$\a^{p}<\a^{p}$, $T^{q}\raj\a^{p}+1 =T^{p}$ etc. Then $P$ is
$\o_{1}$-closed and by CH, it has $\o_{2}$-cc.
$P$ also adds the required tree $T$. (There are no other branches
than $\n_{i}=\cup_{p\in G}\n^{p}_{i}$, $i<\aleph^{\o}$,
because if $p$ forces that
$b$ is a linearly ordered subtree of $T$ and for all $i<\aleph_{\o}$,
$b\ne\n_{i}$, then by using the fact that $P$ is $\o_{1}$-closed, we can find
$q>p$ such that
$q$ decides $b\raj\a^{q}$, $b\raj\a^{q}\ne b^{q}_{i}\raj\a^{q}$
for all $i\in u^{q}$ and no $t\in T^{q}$ is on top of
$b\raj\a^{q}$. Then $q$ forces that the height of $b$ is
$\le\a^{q}$.)
Now, if we want, we can add
$>\aleph_{\o}$ Cohen reals. As in the proof of Claim 4.1,
we can see that this does not add branches to $T$.
$\eop$

\chapter{References}

\item{[Ro]} J. Roitman, The number of automorphism of an atomic
Boolean algebra, Pacific Journal of Mathematics, vol. 94,
1981, 231-242.

\item{[Mo]} J. D. Monk, Automorphism groups, in: J. D. Monk
and R. Bonnet (ed.) Handbook of Boolean Algebras, vol. 2,
North-Holland, Amsterdam, 1989, 517-546.

\bigskip

\settabs\+\hskip 6truecm&\cr

\+Department of mathematics\cr

\+P.O. Box 4\cr

\+00014 University of Helsinki\cr

\+Finland\cr

\bigskip

\+Institute of Mathematics&Rutgers University\cr

\+The Hebrew University&Department of Mathematics\cr

\+Jerusalem&New Brunswick, NJ\cr

\+Israel&USA\cr

\end